\documentclass[11pt]{article}
\usepackage{amsmath} 
\usepackage{amssymb}
\usepackage{amsthm}
\usepackage[usenames]{color}
\usepackage{amscd}
\usepackage{euler}

\definecolor{webgreen}{rgb}{0,.5,0}
\definecolor{webbrown}{rgb}{.6,0,0}

\def\C{{\Bbb C}}

\def\lcm{\operatorname{lcm}}

\begin{document}

\title{\bf A survey of the alternating sum-of-divisors function}
\author{L\'aszl\'o T\'oth \thanks{The author gratefully acknowledges support from the
Austrian Science Fund (FWF) under the projects Nr. P20847-N18 and M1376-N18.}}
\date{}
\maketitle

\centerline{Acta Univ. Sapientiae, Math. {\bf 5} (2013), no. 1, 93--107}

\medskip

\centerline{Dedicated to the memory of my friend and colleague, Professor Antal Bege (1962--2012)}

\begin{abstract} We survey arithmetic and asymptotic properties of the
alternating sum-of-divisors function $\beta$ defined by
$\beta(p^a)=p^a-p^{a-1}+p^{a-2}-\ldots +(-1)^a$ for every prime
power $p^a$ ($a\ge 1$), and extended by multiplicativity. Certain
open problems are also stated.
\end{abstract}

{\it 2010 Mathematics Subject Classification}: 11A25, 11N37

{\it Key words and phrases}: sum-of-divisors function, Liouville function, Euler's totient function,
specially multiplicative function, imperfect number

\section{Introduction}
Let $\beta$ denote the multiplicative arithmetic function defined by
$\beta(1)=1$ and
\begin{equation} \label{def_beta}
\beta(p^a)=p^a-p^{a-1}+p^{a-2}-\ldots +(-1)^a
\end{equation}
for every prime power $p^a$ ($a\ge 1$). That is,
\begin{equation} \label{beta_n}
\beta(n)=\sum_{d\mid n} d\, \lambda(n/d)
\end{equation}
for every integer $n\ge 1$, where $\lambda(n)=(-1)^{\Omega(n)}$ is
the Liouville function, $\Omega(n)$ denoting the number of prime
power divisors of $n$.

The function $\beta$, as a variation of the sum-of-divisors function
$\sigma$, was considered by Martin \cite{Mar1999}, Iannucci
\cite{Ian2006}, Zhou and Zhu \cite{ZhoZhu2009} regarding the
following problem. In analogy with the perfect numbers, $n$ is said
to be imperfect if $2\beta(n) = n$. More generally, $n$ is said to
be $k$-imperfect if $k\beta(n) = n$ for some integer $k\ge 2$. The
only known imperfect numbers are $2,12,40,252,880$, $10\, 880, 75\,
852$, $715\, 816\, 960$ and $3\, 074\, 457\, 344\, 902\, 430\, 720$
(sequence A127725 in \cite{OEIS}). Examples of $3$-imperfect numbers
are $6,120,126,2520$. No $k$-imperfect numbers are known for $k>3$.
See also the book of Guy \cite[p.\ 72]{Guy2004}.

This function occurs in the literature also in another context. Let
\begin{equation*} b(n)=\# \{k: 1\le k\le n \text{ and } \gcd(k,n) \text{ is a square}\}.
\end{equation*}

Then $b(n)=\beta(n)$ ($n\ge 1$), see Cohen \cite[Cor.\
4.2]{Coh1958e1}, Sivaramakrishnan \cite{Siv1979}, \cite[p.\
201]{Siv1989}, McCarthy \cite[Sect.\ 6]{McC1960i}, \cite[p.\
25]{McC1986}, Bege \cite[p.\ 39]{Beg2006}, Iannucci \cite[p.\
12]{Ian2006}. A modality to show the identity $b(n)=\beta(n)$ is to
apply a familiar property of the Liouville function, namely,
\begin{equation} \label{Liouville_prop}
\sum_{d\mid n} \lambda(d) = \chi(n) \quad (n\ge 1),
\end{equation}
where $\chi$ is the characteristic function of the set of squares.
Using \eqref{Liouville_prop},
\begin{equation*}
b(n)= \sum_{k=1}^n \chi(\gcd(k,n))= \sum_{k=1}^n \sum_{d\mid
\gcd(k,n)} \lambda(d)= \sum_{d \mid n} \lambda(d) \sum_{1\le k\le
n,\, d\mid k} 1
\end{equation*}
\begin{equation*}
= \sum_{d\mid n} \lambda(d) \frac{n}{d}=\beta(n).
\end{equation*}

In this paper we survey certain known properties of the function
$\beta$, and give also other ones (without references), which may be
known, but we could not locate them in the literature.

We point out that the function $\beta$ has a double character. On
the one hand, certain properties of this function are similar to
those of the sum-of-divisors function $\sigma$, due to the fact that
both $\beta$ and $\sigma$ are the Dirichlet convolution of two
completely multiplicative functions. Such functions are called in
the literature specially multiplicative functions or quadratic
functions. Their study in connection to the Busche-Ramanujan
identities goes back to the work of Vaidyanathaswamy \cite{Vai1931}.
See also \cite{Hau2003, McC1960ii, McC1986, Siv1989}.

On the other hand, further properties of this function are analogous
to those of the Euler's totient function $\varphi$, as a consequence of the
representation of $\beta$ given above.

We call the function $\beta$ the {\sl alternating sum-of-divisors
function} or {\sl alternating sigma function}. Sivaramakrishnan
\cite[p.\ 201]{Siv1989} remarked that it may be termed the square
totient function.

It is possible, of course, to define other alternating sums of the
positive divisors of $n$. For example, let $\theta(n)=\sum_{d\mid n}
d\, \lambda(d)$ ($n\ge 1$). Then $\theta(n)=\lambda(n)\beta(n)$
($n\ge 1$). This is sequence A061020 in \cite{OEIS}. Another
example: let $n=d_1>d_2>\ldots >d_{\tau(n)}=1$ be the divisors of
$n$, in decreasing order, and let $A(n)=\sum_{j=1}^{\tau(n)}
(-1)^{j-1}d_j$, cf. \cite{Ata2009}. Note that the function $A$ is
not multiplicative.

We do not give detailed proofs, excepting the proofs of formulae
\eqref{prod_id}, \eqref{Rama}, \eqref{asympt_1_beta} and of the
Proposition in Section \ref{Section_Super-imperfect numbers}, which
are included in Section \ref{Section_Proofs}. We leave to the
interested reader to compare the corresponding properties of the
functions $\beta$, $\sigma$ and $\varphi$. See, for example, the
books \cite{Apo1976,HarWri1975,McC1986,Siv1989, Ten1995}.

In Section \ref{Section_Super-imperfect numbers} we pose certain
open problems.  One of them is concerning super-imperfect numbers
$n$, defined by $2\beta(\beta(n))=n$. This notion seems not to
appear in the literature. The super-imperfect numbers up to $10^7$
are $n=2,4,8,128, 32\, 768$. The number $2\, 147\, 483\, 648$ is
also super-imperfect.

The corresponding concept for the sigma function is the following:
A number $n$ is called superperfect if $\sigma(\sigma(n))=2n$. The
even superperfect numbers are $2^{p-1}$, where $2^p-1$ is a Mersenne
prime, cf. \cite{Sur1969} (sequence A019279 in \cite{OEIS}). No odd
superperfect numbers are known.


\section{Basic properties} \label{Section_Basic_properties}
It is clear from \eqref{def_beta} that for every prime power $p^a$
($a\ge 1$),
\begin{equation} \label{beta_pp_formula}
\beta(p^a)= \frac{p^{a+1}+(-1)^a}{p+1} = \begin{cases}
\frac{p^{a+1}-1}{p+1}, & \text{if $a\ge 1$ is odd},\\
\frac{p^{a+1}+1}{p+1}, & \text{if $a\ge 2$ is even}.
\end{cases}
\end{equation}

We obtain from \eqref{beta_n},
\begin{equation*} 
\sum_{n=1}^{\infty} \frac{\beta(n)}{n^s} =
\frac{\zeta(s-1)\zeta(2s)}{\zeta(s)} \quad (\Re (s)>2),
\end{equation*}
leading to another convolution representation of $\beta$, namely
\begin{equation} \label{beta_n_repr}
\beta(n)=\sum_{d^2k=n} \varphi(k) \quad (n\ge 1),
\end{equation}
cf. McCarthy \cite[Sect.\ 6]{McC1960i}, \cite[p.\ 25]{McC1986}, Bege
\cite[p.\ 39]{Beg2006}.

We have $\varphi(n)\le \beta(n)\le n$ ($n\ge 1$). More exactly, it
follows from \eqref{beta_n_repr} that for every $n\ge 1$,
\begin{equation*}  
\beta(n)=\varphi(n)+ \sum_{d^2k=n, d>1} \varphi(k)\ge \varphi(n),
\end{equation*}
with equality for the squarefree values of $n$. Also,
\begin{equation*} 
\beta(n)\le \sum_{dk=n} \varphi(k)=n,
\end{equation*}
with equality only for $n=1$.

Moreover, $\beta(n)\le \varphi^*(n)$ for every $n\ge 1$, with
equality if and only if $n$ is squarefree or twice a squarefree
number. This follows easily from \eqref{beta_pp_formula}. Here
$\varphi^*$ is the unitary Euler function, which is multiplicative
and given by $\varphi^*(p^a)=p^a-1$ for every prime power $p^a$
($a\ge 1$), cf. \cite{McC1986, Siv1989}. Also, $\beta(n)\ge
\sqrt{n}$ ($n\ge 1$, $n\ne 2$, $n\ne 6$).

Similar to the corresponding property of the function $\sigma$,
$\beta(n)$ is odd if and only if $n$ is a square or twice a square.

The function $\beta$ appears in certain elementary identities
regarding the set of squares, for example in
\begin{equation*}
\sum_{\substack{k=1\\ \gcd(k,n) \text{ a square}}}^n k =
\frac{n(\beta(n)+\chi(n))}{2} \quad (n\ge 1),
\end{equation*}
\begin{equation*}
\prod_{\substack{k=1\\ \gcd(k,n) \text{ a square}}}^n k =
n^{\beta(n)} \prod_{d\mid n} (d!/d^d)^{\lambda(n/d)} \quad (n\ge 1),
\end{equation*}
which can be deduced from \eqref{Liouville_prop}.


\section{Generalizations} \label{Section_Generalizations}

An obvious generalization of $\beta$ is the function $\beta_a$
($a\in \C$) defined by
\begin{equation} \label{beta_a}
\beta_a(n)=\sum_{d\mid n} d^a \, \lambda(n/d) \quad (n\ge 1).
\end{equation}

If $a=m$ is a positive integer, then the following representation
can be given: $\beta_m(n)= \# \{k: 1\le k\le n^m, (k,n^m)_m \text{
is a $2m$-th power}\}$, where $(a,b)_m$ stands for the largest
common $m$-th power divisor of $a$ and $b$. See McCarthy
\cite[Sect.\ 6]{McC1960i}, \cite[p.\ 51]{McC1986}.

Note that if $a=0$, then $\beta_0=\chi$, the characteristic function
of the set of squares, used above.

For an arbitrary nonempty set $S$ of positive integers let
$\varphi_S(n)= \# \{k: 1\le k\le n, \gcd(k,n)\in S\}$. For $S=\{1\}$
and $S$ the set of squares this reduces to Euler's function
$\varphi$ and to the function $\beta$, respectively. The function
$\varphi_S$ was investigated by Cohen \cite{Coh1959}. For every set
$S$ one has
\begin{equation*}
\varphi_S(n)= \sum_{d\mid n} d\mu_S(n/d) \quad (n\ge 1),
\end{equation*}
where the function $\mu_S$ is defined by $\sum_{d\mid n} \mu_S(d)
=\chi_S(n)$ ($n\ge 1$), i.e., $\mu_S=\mu*\chi_S$ in terms of the
Dirichlet convolution $*$, $\chi_S$ and $\mu$ denoting the
characteristic function of $S$ and the M\"obius function,
respectively.

Also, let
\begin{equation*}
B(r,n)=\sum_{\substack{k=1\\ \gcd(k,n) \text{ a square}}}^n
\exp(2\pi ikr/n),
\end{equation*}
which is an analog of the Ramanujan sum to be considered in
Section \ref{Section_Asymptotic behavior}. Then
\begin{equation*}
B(r,n)= \sum_{d\mid \gcd(r,n)} d\, \lambda(n/d) \quad (r,n\ge 1),
\end{equation*}
see Sivaramakrishnan \cite{Siv1979}, \cite[p.\ 202]{Siv1989}, Haukkanen \cite{Hau1997}. For
$r=n$ one has $B(n,n)=\beta(n)$.

These generalizations can also be combined. See also Haukkanen
\cite{Hau1989,Hau1995}. Many of the results given in the present
paper can be extended for these generalizations.

We consider in what follows only the functions $\beta_a$ defined by
\eqref{beta_a} and do not deal with other generalizations.


\section{Further properties} \label{Further_properties}

For every $n,m\ge 1$,
\begin{equation} \label{Busche_Raman_i}
\beta(n)\beta(m) = \sum_{d\mid \gcd(n,m)} \beta(nm/d^2)d\,
\lambda(d),
\end{equation}
and equivalently,
\begin{equation} \label{Busche_Raman_ii}
\beta(nm) = \sum_{d\mid \gcd(n,m)} \beta(n/d)\beta(m/d)d\mu^2(d),
\end{equation}
cf. \cite{Siv1979}, \cite[p.\ 26]{McC1986}. Here
\eqref{Busche_Raman_i} and \eqref{Busche_Raman_ii} are special cases
of the Busche-Ramanujan identities, valid for specially
multiplicative functions. See \cite{Hau2003,McC1960ii,McC1986,
Siv1989,Vai1931} for their discussions and proofs.

Direct proofs of \eqref{Busche_Raman_i} and \eqref{Busche_Raman_ii}
can be given by showing that both sides of these identities are
multiplicative, viewed as functions of two variables and then
computing their values for prime powers. Recall that an arithmetic
function $f$ of two variables is called multiplicative if it is
nonzero and $f(n_1m_1,n_2m_2)= f(n_1,n_2) f(m_1,m_2)$ holds for any
$n_1,n_2,m_1,m_2\ge 1$ such that $\gcd(n_1n_2,m_1m_2)=1$. See
\cite{Vai1931}, \cite{Tot2011}, \cite[Ch.\ VII]{Siv1989}.

The proof of the equivalence of identities of type
\eqref{Busche_Raman_i} and \eqref{Busche_Raman_ii} is outlined in
\cite{Hau1999}, referring to the work of Vaidyanathaswamy
\cite{Vai1931}.

It follows at once from \eqref{Busche_Raman_ii} that $\beta(nm)\ge
\beta(n)\beta(m)$ for every $n,m\ge 1$, i.e., $\beta$ is
super-multiplicative. Formula \eqref{Busche_Raman_ii} leads also to
the double Dirichlet series
\begin{equation*}
\label{double_Dirichlet} \sum_{n,m=1}^{\infty}
\frac{\beta(nm)}{n^sm^t} =
\frac{\zeta(s-1)\zeta(2s)\zeta(t-1)\zeta(2t)\zeta(s+t-1)}{\zeta(s)\zeta(t)\zeta(2(s+t-1))},
\end{equation*}
valid for $s,t\in \C$ with $\Re (s)>2, \Re (t)>2$.

The generating power series of $\beta$ is
\begin{equation*}
\sum_{n=1}^{\infty} \beta(n) x^n = \sum_{n=1}^{\infty}
\frac{\lambda(n)x^n}{(1-x^n)^2} \quad (|x|<1),
\end{equation*}
which is a direct consequence of \eqref{beta_n}.

Consider the functions $\beta_a$ defined by \eqref{beta_a}. One has
\begin{equation} \label{beta_a_beta_b}
\sum_{n=1}^{\infty} \frac{\beta_a(n)\beta_b(n)}{n^s} =
\frac{\zeta(s)\zeta(s-a-b)\zeta(2s-2a)\zeta(2s-2b)}{\zeta(s-a)\zeta(s-b)\zeta(2s-a-b)},
\end{equation}
valid for every $s,a,b\in \C$ with $\Re (s)>1+\max(0,\Re (a), \Re
(b), \Re (a+b))$.

This formula is similar to Ramanujan's well-known result for the
product $\sigma_a(n)\sigma_b(n)$, where $\sigma_a(n)=\sum_{d\mid n}
d^a$. Formula \eqref{beta_a_beta_b} is due to Chowla \cite{Cho1928},
in an equivalent form for the product $\theta_a(n)\theta_b(n)$,
where $\theta_a(n)=\sum_{d\mid n} d^a\, \lambda(d)$.

Formula \eqref{beta_a_beta_b} and that of Ramanujan follow from the
next more general result concerning the product of two arbitrary
specially multiplicative functions.

If $f,g,h,k$ are completely multiplicative functions, then
\begin{equation} \label{prod_id}
(f*g)(h*k)=fh*fk*gh*gk*w,
\end{equation}
where $w(n)=\mu(m)f(m)g(m)h(m)k(m)$ if $n=m^2$ is a square and
$w(n)=0$ otherwise.

This result is given by Vaidyanathaswamy \cite[p.\ 621]{Vai1931}, Lambek
\cite{Lam1966}, Subbarao \cite{Sub1968}. See also \cite[p.\
50]{Siv1989}. The proof of \eqref{prod_id} can be carried out using
Euler products. This is well known in the case of Ramanujan's result
regarding $\sigma_a \sigma_b$, and is presented in several texts, cf.,
e.g., \cite[Th.\ 305]{HarWri1975}, \cite[Prop.\ 5.4]{McC1986}. An
alternative proof is given by Lambek \cite{Lam1966}.

In Section \ref{Section_Proofs} we offer another less known short
proof of \eqref{prod_id}.

In the case of the functions $f(n)=n^a$, $h(n)=n^b$,
$g(n)=k(n)=\lambda(n)$ we obtain \eqref{beta_a_beta_b} by using the
known formulae for the Dirichlet series corresponding to the right
hand side of \eqref{prod_id}.

If $f(n)=n^a$, $h(n)=n^b$, $g(n)=\lambda(n)$, $k(n)=1$, then we
deduce
\begin{equation*}
\sum_{n=1}^{\infty} \frac{\beta_a(n)\sigma_b(n)}{n^s} =
\frac{\zeta(s-a)\zeta(s-a-b)\zeta(2s)\zeta(2s-2b)\zeta(2s-a-b)}{\zeta(s)\zeta(s-b)\zeta(4s-2a-2b)},
\end{equation*}
valid for the same region as \eqref{beta_a_beta_b}.

Remark that we obtain, as direct corollaries, the next formulae:
\begin{equation*}
\sum_{n=1}^{\infty} \frac{\beta^2(n)}{n^s} =
\frac{\zeta(s)\zeta(s-2)\zeta(2s-2)}{\zeta^2(s-1)},
\end{equation*}
\begin{equation} \label{beta_n_2}
\sum_{n=1}^{\infty} \frac{\beta(n^2)}{n^s} =
\frac{\zeta(s)\zeta(s-2)}{\zeta(s-1)},
\end{equation}
\begin{equation*}
\sum_{n=1}^{\infty} \frac{\beta(n)\sigma(n)}{n^s} =
\frac{\zeta(s-2)\zeta(2s)\zeta^2(2s-2)}{\zeta(s)\zeta(4s-4)},
\end{equation*}
all valid for $\Re (s)>3$. Here \eqref{beta_n_2} is obtained from
\eqref{beta_a_beta_b} by choosing $a=1$ and $b=0$.

From these Dirichlet series representations we can deduce the
following convolutional identities:
\begin{equation} \label{convo_beta_2}
\beta^2(n)= \sum_{dk=n} d\, 2^{\omega(d)}\lambda(d) \sigma_2(k)
\quad (n\ge 1),
\end{equation}
\begin{equation} \label{convo_beta_n_2}
\beta(n^2) = \sum_{dk=n} d \mu(d) \sigma_2(k) \quad (n\ge 1),
\end{equation}
\begin{equation} \label{convo_beta_sigma}
\beta(n)\sigma(n)= \sum_{d^2k=n} d^2 2^{\omega(d)}\beta_2(k) \quad
(n\ge 1),
\end{equation}
where $\omega(n)$ denotes the number of distinct prime factors of
$n$.


\section{Asymptotic behavior} \label{Section_Asymptotic behavior}

The average order of $\beta(n)$ is $(\pi^2/15)n$, more exactly,
\begin{equation} \label{asympt_beta}
\sum_{n\le x} \beta(n) = \frac{\pi^2}{30}x^2 + {\cal O} \left(x(\log
x)^{2/3}(\log \log x)^{4/3} \right).
\end{equation}

Formula \eqref{asympt_beta} follows from the convolution
representation \eqref{beta_n_repr} and from the known estimate of
Walfisz regarding $\sum_{n\le x} \varphi(n)$ with the same error
term as above.

There are also other asymptotic properties of the $\varphi$
function, which can be transposed to $\beta$ by using that
$\beta(n)\ge \varphi(n)$, with equality for $n$ squarefree.  For
example,
\begin{equation*}
\liminf_{n\to \infty} \frac{\beta(n)\log \log n}{n} =e^{-\gamma},
\end{equation*}
where $\gamma$ is Euler's constant (cf. \cite[Th.\ 328]{HarWri1975}
concerning $\varphi$). Another example: the set $\{\beta(n)/n: n\ge
1 \}$ is dense in the interval $[0,1]$.

Let $c_r(n)$ denote the Ramanujan sum, defined as the sum of $n$-th
powers of the primitive $r$-th roots of unity. Then
\begin{equation} \label{Rama}
\frac{\beta(n)}{n}= \frac{\pi^2}{15} \sum_{r=1}^{\infty}
\frac{\lambda(r)}{r^2}c_r(n)
\end{equation}
\begin{equation*}
= \frac{\pi^2}{15} \left(1-\frac{(-1)^n}{2^2}- \frac{2\cos(2\pi
n/3)}{3^2} + \frac{2\cos(\pi n/2)}{4^2}+\ldots \right),
\end{equation*}
showing how the values of $\beta(n)/n$ fluctuate harmonically about
their mean value $\pi^2/15$, cf. \cite{Coh1961}, \cite[p.\
245]{McC1986}.

A quick direct proof of formula \eqref{Rama} is given in Section
\ref{Section_Proofs}. We refer to \cite{Luc2010} for a recent survey
of expansions of functions with respect to Ramanujan sums.

From the identities \eqref{convo_beta_2}, \eqref{convo_beta_n_2} and
\eqref{convo_beta_sigma} we deduce the following asymptotic
formulae:
\begin{equation*}
\sum_{n\le x} \beta^2(n) = \frac{2\zeta(3)}{15}x^3 + {\cal O}
\left(x^2 (\log x)^2 \right),
\end{equation*}
\begin{equation*}
\sum_{n\le x} \beta(n^2) = \frac{2\zeta(3)}{\pi^2}x^3 + {\cal O}
\left(x^2 \log x \right),
\end{equation*}
\begin{equation*}
\sum_{n\le x} \beta(n)\sigma(n) = \frac{\pi^6}{2430\zeta(3)}x^3 +
{\cal O} \left(x^2 \right).
\end{equation*}

We also have
\begin{equation} \label{asympt_1_beta}
\sum_{n\le x} \frac1{\beta(n)} = K_1\log x + K_2 + {\cal O}
\left(x^{-1+\varepsilon}\right),
\end{equation}
for every $\varepsilon >0$, where $K_1$ and $K_2$ are constants,
\begin{equation*}
K_1= \prod_p \left(1-\frac1{p}\right)\left(1+\sum_{a=1}^{\infty}
\frac1{\beta(p^a)}\right).
\end{equation*}

For the proof of \eqref{asympt_1_beta} see Section
\ref{Section_Proofs}.


\section{Unitary analog} \label{Section_Unitary_analog}

Consider the function $\beta^*$ defined by
\begin{equation*}
\beta^*(n)=\sum_{d\mid\mid n} d\, \lambda(n/d) \quad (n\ge 1),
\end{equation*}
where the sum is over the unitary divisors $d$ of $n$. Recall that
$d$ is a unitary divisor of $n$ if $d\mid n$ and $\gcd(d,n/d)=1$.
Here $\beta^*(p^a)=p^a+(-1)^a$ for every prime power $p^a$ ($a\ge
1$) and
\begin{equation} \label{series_beta_*}
\sum_{n=1}^{\infty} \frac{\beta^*(n)}{n^s} =
\frac{\zeta(s-1)\zeta(2s)\zeta(2s-1)}{\zeta(s)\zeta(4s-2)} \quad
(\Re (s) > 2).
\end{equation}

The formula \eqref{series_beta_*} can be derived using Euler
products or by establishing the convolutional identity
\begin{equation*} 
\beta^*(n)= \sum_{dk^2=n} \beta(d) k\, q(k),
\end{equation*}
$q$ standing for the characteristic function of the squarefree
numbers. This leads also to the asymptotic formula
\begin{equation*}
\sum_{n\le x} \beta^*(n) = \frac{63\zeta(3)}{2\pi^4}x^2 + {\cal O}
\left(x(\log x)^{5/3}(\log \log x)^{4/3} \right).
\end{equation*}

Note the following interpretation: $\beta^*(n)=\# \{k: 1\le k\le n
\text{ and } \gcd(k,n)_*$ $\text{is a square}\}$, where $(a,b)_*$ is
the largest divisor of $a$ which is a unitary divisor of $b$.


\section{Super-imperfect numbers, open problems} \label{Section_Super-imperfect numbers}

A number $n$ is super-imperfect if $2\beta(\beta(n))=n$, cf.
Introduction. Observe that, excepting $4$, all the other examples of
super-imperfect numbers are of the form $n=2^{2^k-1}$ with $k\in
\{1,2,3,4,5\}$. The proof of the next statement is given in Section
\ref{Section_Proofs}.

\textit{Proposition}. For $k\ge 1$ the number $n_k=2^{2^k-1}$ is
super-imperfect if and only if $k\in \{1,2,3,4,5\}$.

\textit{Problem} 1. Is there any other super-imperfect number?

More generally, we define $n$ to be $(m,k)$-imperfect if
$k\beta^{(m)}(n)=n$, where $\beta^{(m)}$ is the $m$-fold iterate of
$\beta$. For example, $3$, $15$, $18$, $36$, $72$, $255$ are
$(2,3)$-imperfect, $6$, $12$, $24$, $30$, $60$, $120$
are $(2,6)$-imperfect, $6$, $36$, $144$ are $(3,6)$-imperfect numbers.

We refer to \cite{CohRie1996} regarding $(m,k)$-perfect numbers,
defined by $\sigma^{(m)}(n)=kn$.

\textit{Problem} 2. Investigate the $(m,k)$-imperfect numbers.

The numbers $n=1,20,116,135,171,194,740,\ldots$ are solutions of the equation $\beta(n)=\beta(n+1)$.

\textit{Problem} 3. Are there infinitely many numbers $n$ such that
$\beta(n)=\beta(n+1)$?

Remark that it is not known if there are infinitely many numbers $n$
such that $\sigma(n)=\sigma(n+1)$ (sequence A002961 in \cite{OEIS}).
See also Weingartner \cite{Wei2011}.

The next problem is the analog of Lehmer's open problem concerning the $\varphi$ function.

\textit{Problem} 4. Is there any composite number $n\ne 4$ such that
$\beta(n)$ divides $n-1$?

Up to $10^6$ there are no such composite numbers.

The computations were performed using Maple. The function $\beta(n)$
was generated by the following procedure:
\begin{verbatim}
  beta:= proc(n) local x, i: x:= 1:
  for i from 1 to nops(ifactors(n)[2])
  do p_i:= ifactors(n)[2][i][1]: a_i:= ifactors(n)[2][i][2];
  x:= x*((p_i^(a_i+1)+(-1)^(a_i))/(p_i+1)): od: RETURN(x) end;
  # alternating sum-of-divisors function
\end{verbatim}


\section{Proofs} \label{Section_Proofs}

\textit{Proof of formula} \eqref{prod_id}: Write
\begin{equation*}
(f*g)(n)(h*k)(n)= \sum_{\substack{d\mid n\\e\mid n}}
f(d)g(n/d)h(e)k(n/e),
\end{equation*}
where $d\mid n, e\mid n \ \Leftrightarrow \ \lcm[d,e] \mid n$. Write
$d=mu$, $e=mv$ with $\gcd(u,v)=1$. Then $\lcm[d,e] =muv$ and obtain
that this sum is
\begin{equation*}
\sum_{\substack{muv\mid n \\ \gcd(u,v)=1}}
f(mu)g(n/(mu))h(mv)k(n/(mv))
\end{equation*}
\begin{equation*}
= \sum_{muv\mid n} f(mu)g(n/(mu))h(mv)k(n/(mv)) \sum_{\delta \mid
(\gcd(u,v)} \mu(\delta).
\end{equation*}

Putting now $u=\delta x$, $v=\delta y$ and using that the considered
functions are all completely multiplicative the latter sum is
\begin{equation*}
\sum_{\delta^2xymt=n} (\mu fghk)(\delta) (fk)(x) (gh)(y) (fh)(m)
(gk)(t),
\end{equation*}
finishing the proof (cf. \cite[p.\ 621]{Vai1931} and \cite[p.\
161]{Ten1995}).

\textit{Proof of formula} \eqref{Rama}: Let $\eta_r(n)=r$ if $r\mid
n$ and $\eta_r(n)=0$ otherwise. Applying that $\sum_{d\mid r}
c_d(n)=\eta_r(n)$ we deduce
\begin{equation*}
\frac{\beta(n)}{n}=\sum_{d\mid n} \frac{\lambda(d)}{d}=
\sum_{d=1}^{\infty} \frac{\lambda(d)}{d^2}\eta_d(n)=
\sum_{d=1}^{\infty} \frac{\lambda(d)}{d^2} \sum_{r\mid d} c_r(n)
\end{equation*}
\begin{equation*}
= \sum_{r=1}^{\infty} \frac{\lambda(r)}{r^2} c_r(n)
\sum_{k=1}^{\infty} \frac{\lambda(k)}{k^2} =
\frac{\zeta(4)}{\zeta(2)} \sum_{r=1}^{\infty} \frac{\lambda(r)}{r^2}
c_r(n),
\end{equation*}
using that $\lambda$ is completely multiplicative and its Dirichlet
series is $\sum_{n=1}^{\infty} \lambda(n)/n^s$ $=
\zeta(2s)/\zeta(s)$. The rearranging of the terms is justified by
the absolute convergence.

\textit{Proof of formula} \eqref{asympt_1_beta}: Write
\begin{equation*}
\frac1{\beta(n)} = \sum_{\substack{dk=n\\ \gcd(d,k)=1}}
\frac{h(d)}{\varphi(k)}
\end{equation*}
as the unitary convolution of the functions $h$ and $1/\varphi$,
where $h$ is multiplicative and for every prime power $p^a$ ($a \ge
1$),
\begin{equation*}
\frac1{\beta(p^a)} = h(p^a)+ \frac1{\varphi(p^a)}, \quad h(p^a)= -
\frac{p^{a-1}+(-1)^a}{p^{a-1}(p-1)(p^{a+1}+(-1)^a)}.
\end{equation*}

Here
\begin{equation*}
|h(p^a)|\le \frac1{p^a(p-1)^2}, \quad |h(n)|\le
\frac{f(n)}{\varphi(n)} \ (n\ge 1),
\end{equation*}
with $f(n)=\prod_{p\mid n} (p(p-1))^{-1}$. We deduce
\begin{equation*}
\sum_{n\le x} \frac1{\beta(n)} = \sum_{d\le x} h(d)
\sum_{\substack{k\le x/d \\ \gcd(d,k)=1}} \frac1{\varphi(k)},
\end{equation*}
and use the known estimates for the inner sum. The same arguments
were applied in the proof of \cite[Th.\ 2]{Tot2010}.

\textit{Proof of the Proposition of Section}
\ref{Section_Super-imperfect numbers}: The fact that the numbers
$n_k$ with $1\le k\le 5$ are super-imperfect follow also by direct
computations, but the following arguments reveal a connection to the
Fermat numbers $F_m=2^{2^m}+1$.

For $n_k=2^{2^k-1}$ with $k\ge 1$ we have
\begin{equation*}
\beta(n_k)= \frac{2^{2^k}-1}{3} = F_1F_2\cdots F_{k-1}
\end{equation*}
(for $k=1$ this is $1$, the empty product). Since the numbers $F_m$
are pairwise relatively prime,
\begin{equation*}
\beta(\beta(n_k))= \beta(F_1)\beta(F_2)\cdots \beta(F_{k-1}).
\end{equation*}

Now for $2\le k\le 5$, using that $F_1, F_2, F_3,F_4$ are primes,
\begin{equation*}
\beta(\beta(n_k))= 2^{2^1}\cdot 2^{2^2}\cdots 2^{2^{k-1}}= 2^{2^k-2}
= \frac{n_k}{2},
\end{equation*}
showing that $n_k$ is super-imperfect.

Now let $k\ge 6$. We use that $F_5$ is composite and that $\beta(n)
\lneqq n-1$ for every $n\ne 4$ composite. Hence $\beta(F_5)\lneqq
2^{2^5}$ and
\begin{equation*}
\beta(\beta(n_k)) \lneqq \beta(F_1)\beta(F_2)\cdots \beta(F_{k-1})=
2^{2^k-2}=\frac{n_k}{2},
\end{equation*}
ending the proof.

Note that for $k\ge 2$ the number $m_k= 2^{2^k-1} F_1F_2\cdots
F_{k-1}$ is imperfect if and only if $k\in \{2,3,4,5\}$.  This
follows by similar arguments. The imperfect numbers of this form are
$40$, $10\, 880$, $715\, 816\, 960$ and $3\, 074\, 457\, 344\, 902\,
430\, 720$, given in the Introduction.


\noindent L. T\'oth \\
Institute of Mathematics, Universit\"at f\"ur Bodenkultur \\
Gregor Mendel-Stra{\ss}e 33, A-1180 Vienna, Austria \\ and \\
Department of Mathematics, University of P\'ecs \\ Ifj\'us\'ag u. 6,
H-7624 P\'ecs, Hungary \\ E-mail: ltoth@gamma.ttk.pte.hu

\end{document}